%Paper: funct-an/9506003
%From: guido@mat.utovrm.it (Daniele Guido, U. Roma Tor Vergata)
%Date: Wed, 14 Jun 1995 19:42:00 +0100

%%%%%% FORMAT %%%%%%%%%%%%

\magnification\magstep1
\hoffset=0.5truecm
\voffset=0.5truecm
\hsize=15.8truecm
\vsize=23.truecm
\baselineskip=14pt plus0.1pt minus0.1pt \parindent=19pt
\lineskip=4pt\lineskiplimit=0.1pt      \parskip=0.1pt plus1pt

%%%%%%%% Macros for theorems %%%%%%%%
\font\sectionfont=cmbx10 scaled\magstep1
\def\titlea#1{\vskip0pt plus.3\vsize\penalty-75
    \vskip0pt plus -.3\vsize\bigskip\bigskip
    \noindent{\sectionfont #1}\nobreak\smallskip\noindent}

\def\claim#1#2{\vskip.1in\medbreak\noindent{\bf #1.} {\sl #2}\par
    \ifdim\lastskip<\medskipamount\removelastskip\penalty55\medskip\fi}
\def\beglemma#1#2\endlemma{\claim{#1 Lemma}{#2}}
\def\begdefinition#1#2\enddefinition{\claim{#1 Definition}{#2}}
\def\begtheorem#1#2\endtheorem{\claim{#1 Theorem}{#2}}
\def\begcorollary#1#2\endcorollary{\claim{#1 Corollary}{#2}}
\def\begremark#1#2\endremark{\vskip.1in\medbreak\noindent{\bf#1 Remark.} #2
    \par\ifdim\lastskip<\medskipamount\removelastskip\penalty55\medskip\fi}
\def\begproposition#1#2\endproposition{\claim{#1 Proposition}{#2}}
\def\begProof{\medskip\noindent{\bf Proof.}\quad}
\def\begProofof#1{\medskip\noindent{\bf Proof of #1.}\quad}
\def\square{\hbox{$\sqcap\!\!\!\!\sqcup$}}
\def\Halmos{\hfill\square}
\def\endProof{\Halmos\par
    \ifdim\lastskip<\medskipamount\removelastskip\penalty55\medskip\fi}
\newcount\FNOTcount \FNOTcount=1

\def\addfnot{\global\advance\FNOTcount by 1}

%%%% Reference macroes %%%%%%%%%%%%%%%%%%%%
\font\sectionfont=cmbx10 scaled\magstep1
\def\titlea#1{\vskip0pt plus.3\vsize\penalty-75
    \vskip0pt plus -.3\vsize\bigskip\bigskip
    \noindent{\sectionfont #1}\nobreak\smallskip\noindent}
\def\refno#1#2{\item{[#1]}{#2}}
\def\begref#1#2{\titlea{#1}}
\def\endref{}
\def\np{\par\noindent}

\newcount\REFcount \REFcount=1
\def\numref{\number\REFcount}
\def\addref{\global\advance\REFcount by 1}
\def\wdef#1#2{\expandafter\xdef\csname#1\endcsname{#2}}
\def\wdch#1#2#3{\ifundef{#1#2}\wdef{#1#2}{#3}
    \else\write16{!!doubly defined#1,#2}\fi}
\def\wval#1{\csname#1\endcsname}
\def\ifundef#1{\expandafter\ifx\csname#1\endcsname\relax}

\def\autonumref{
    \def\rfr(##1){\wdef{q##1}{yes}\ifundef{r##1}$\diamondsuit$##1
        \write16{!!ref ##1 was never defined!!}\else\wval{r##1}\fi}
    \def\REF(##1)##2\endREF{\wdch{r}{##1}{\numref}\addref}\REFERENCES
    \def\references{
        \def\REF(####1)####2\endREF{
            \ifundef{q####1}\write16{!!ref. [####1] was never quoted!!}\fi
            \refno{\rfr(####1)}####2}
        \begref{References}{99}\REFERENCES\endref}}

\def\REFERENCES{
 \REF(AGPS4) Albeverio S., Guido D., Ponosov A., Scarlatti S., ``{\it Singular
traces and compact operators, I}'', to appear in Journal of Functional
Analysis.
 \endREF
 \REF(Bess1) Besse A.L., ``{\it Manifolds all of whose geodesics are closed}'',
Springer Verlag, Berlin Heidelberg New York, 1978.
 \endREF
 \REF(Brat1) Bratteli O., ``{\it Derivations, Dissipations and Group actions on
$C^*$-algebras.}'' Lecture Notes in Math. 1229. Springer, New York 1986.
 \endREF
 \REF(BrRo1) Bratteli O., Robinson D.W., ``{\it Operator algebras and quantum
statistical mechanics, I.}'' Sprin\-ger, New York, 1979.
 \endREF
 \REF(Conn4) Connes A., ``{\it Non Commutative Geometry}'', Academic Press, New
York, San Francisco, London, 1994.
 \endREF
 \REF(Conn5) Connes A., ``{\it Compact metric spaces, Fredholm modules and
Hyperfinitness}'', Ergod. Theory Dinam. Sys. {\bf 9} (1989), 207-230.
 \endREF
 \REF(CoMo1) Connes A., Moscovici ,``{\it Transgression du caract\`ere de
Chern et cohomologie cyclique}'', C.R. Acad. Sci. Paris Ser. I Math. {\bf 303}
(1986), 913-918.
 \endREF
 \REF(Dixm1) Dixmier J., ``{\it Existence de traces non normales}'', C.R. Acad.
Sci. Paris {\bf 262}, (1966).
 \endREF
 \REF(DuGu1) Duistermaat J.J., Guillemin V.W., ``{\it The spectrum of positive
elliptic operators and periodic bicharacteristics}'' Invent. Math. {\bf 29},
(1975) 39-79.
 \endREF
 \REF(GBVa1) V\'arilly J. C., Gracia-Bond{\'\i}a J. M., ``{\it Connes'
Noncommutative Differential geometry and the Standard Model}", J. Geom. Phys.
{\bf 12}, (1993) 223.
 \endREF
 \REF(GuIs1) Guido D., Isola T., ``{\it Singular traces on semifinite
von~Neumann algebras}'' to appear in Journal of Functional Analysis.
 \endREF
 \REF(Voic1) Voiculescu D., ``{\it On the existence of quasicentral approximate
units relative to normed ideals. Part I}'', Journ. Funct. An. {\bf 91}, (1990)
1-36.
 \endREF
 }
 \autonumref

%%%%%%%%%%%%%%%%%%%%%%%%%%%%%%%%%%%%%%%%%%%%

%%%%%% Special Fonts %%%%%%%
%Number Sets     % Blackboard bold
\def\Co{{\bf C}} %\def\Co{{\Bbb C}}
\def\Na{{\bf N}} %\def\Na{{\Bbb N}}
\def\Re{{\bf R}} %\def\Re{{\Bbb R}}
 %\def\To{{\Bbb T}}
 %\def\Ze{{\Bbb Z}}

%%%%%% SYMBOLS %%%%%%%%%%%%%
\def\A{{\cal A}}
\def\a{\alpha}
\def\B{{\cal B}}

\def\C{{\cal C}}
\def\d{\delta}
\def\D{{\cal D}}
\def\f{\varphi}

\def\H{{\cal H}}

\def\K{{\cal K}}

\def\l{\lambda}

\def\m{\mu}

\def\o{\omega}

\def\Q{\Omega}

\def\s{\sigma}

\def\t{\tau}

\def\w{\omega}

\def\norm#1{|\kern-.5mm |\kern-.5mm | #1 |\kern-.5mm |\kern-.5mm |}
\def\ov{\overline}

%%% First page %%%

%%% TITLE %%%
\topskip3.cm
\font\ftitle=cmbx10 scaled\magstep1
\centerline{\ftitle A REMARK ON TRACE PROPERTIES OF K-CYCLES.}
 \bigskip\bigskip
\centerline{Fabio Cipriani$^1$\footnote{$^*$}
{Supported by H.C.M. EC-grant}
Daniele Guido$^2$\footnote{$^{**}$}
{ Supported in part by MURST and
CNR-GNAFA.} and Sergio Scarlatti$^{2}$}
\footnote{}{E-mail:\ fabio@maths.nott.ac.uk, guido@mat.utovrm.it,
scarlatti@mat.utovrm.it }
 \vskip1.truecm
\item{$(^1)$} Mathematical Department, University of Nottingham,
\par
University Park NG7-2RD  , Nottingham UK
\item{$(^2)$}
 Dipartimento di Matematica, Universit\`a di Roma ``Tor Vergata''
 \par
 Via della Ricerca Scientifica, I--00133 Roma, Italia.
%%% end TITLE %%%
\vskip2.cm
%\hfill June 1995\vskip2.cm
\noindent
 {\bf Abstract.} In this paper we discuss trace properties of $d^+$-summable
$K$-cycles considered by A.Connes in [\rfr(Conn4)].
 More precisely we give a proof of a trace theorem on the algebra $\A$ of a
$K$--cycle stated in [\rfr(Conn4)], namely we show that a natural functional on
$\A$ is a trace functional. Then we discuss whether this functional gives a
trace on the whole universal graded differential algebra $\Q(\A)$.
On the one hand we prove that the regularity conditions on $K$-cycles
considered
 in [\rfr(Conn4)] imply the trace property on $\Q(\A)$.
On the other hand, by constructing
an explicit counterexample, we remark that the sole $K$-cycle assumption
is not sufficient for such a property to hold.

\vfill\eject

\topskip0.cm

\titlea{Introduction}

A major role in A. Connes' non commutative geometry is played by the concept
of  $d^+$--summable $K$--cycle $(\A,D,\H)$. Such an object
generalizes, in an operator--theoretical framework, classical Riemannian
geometry on a compact space. There, given a smooth, closed, compact, Riemannian
spin--manifold $(M,g)$,  one has the Riemannian measure $dm_g$
over $M$ and the Hilbert space  $\H=L^2(M,S,dm_g)$ of the
square--integrable sections of the spinor bundle $S$ on $M$. On $\H$ there is
a natural action  of the *--algebra $\A$ of C$^{\infty}$--functions (by
pointwise multiplication of sections) and the action of the Dirac operator
$D=\partial_M$, associated with a fixed Riemannian connection over $S$ (in
order
to simplify the discussion, we shall assume $Ker(\partial_M)=\{0\}$).
 A key relation among these objects is that $D$ implements a $B(\H)$--valued
derivation on $\A$, i.e. the commutator $[D,f]$ between $D$ and any element
 $f\in\A$ seen as an operator on $\H$ is a bounded operator (given by the
Clifford multiplication by the external derivative of $f$).
 In this operator--theoretical framework the dimension of the manifold $M$ can
be detected looking at the eigenvalue distribution of the Laplacian operator
$D^2$. In fact the celebrated Weyl's theorem implies $\sup{1 \over \log
n}\sigma_n (|D|^{-d})<\infty$ (where $d$ is the dimension of $M$ and $\sigma_n
(T)$ denotes the sum of the first  $n$ eigenvalues $\lambda_k (|T|)$ of  the
compact operator $|T|$).
 In the  language of the normed ideals this can be expressed saying that
$|D|^{-d}$ belongs to the dual ${\cal L}^{1^+}(\H )$ of the Macaev ideal. The
trace--theorem of A. Connes ([\rfr(Conn5)]) allows to reconstruct, from
the data $(\A,D,\H)$, the  Riemannian measure of $M$:
 $$
\int_M fdm_g =c(d)\cdot  \tau_\o (f |D|^{-d})
 $$
where $\tau_\o$ is the Dixmier trace (a non normal trace over $B(\H)$ which is
finite precisely on ${\cal L}^{1^+}(\H )$; see  below) and $c(d)$ is a constant
depending only on $d$.

 A. Connes stated that for any $K$--cycle $(\A,D,\H)$ the state
$\f(a):=\t_\o (|D|^{-d}a)$, defined on $B(\H)$ gives a hypertrace on $\A$
([\rfr(Conn4)], Proposition 15 Chapter IV.2). He also discussed some
consequences of this result on the hermitian structures and the action
functional associated with the $K$--cycle ([\rfr(Conn4)], Proposition 5 Chapter
VI.1).
 To our knowledge, the proof of the previous statement has not been published.
On the contrary, in some published papers such result appeared as an assumption
(see e.g. [\rfr(GBVa1)]). Therefore we think it is worthwhile to give a proof
of it (theorem 1.2 of the next section).

 The functional $\f$ gives rise to a state on the whole universal graded
differential algebra $\Q(\A)$. In the commutative case this is the trace on
the Clifford bundle of the manifold $M$, hence it is natural to ask whether the
trace property for such a state holds for the algebra $\Q(\A)$ associated with
any $K$--cycle.

 The preceding result does not hold in full generality, indeed
Section~2 is devoted to the construction of a counterexample of
such a statement.
 However we shall prove the trace property on $\Q(\A)$ under
natural regularity conditions (Theorem~1.6).
 These conditions have been proposed by Connes [\rfr(Conn4),
p.546] (cf. also [\rfr(Conn4), Theorem~8 p.308]) as the analogue of
differentiability conditions for functions and differential forms
in the non commutative setting.

\titlea{1. Traces constructed by K-cycles}

 We start recalling the definition of $K$-cycle given by A. Connes
[\rfr(Conn4)]. Let $B(\H)$ be the space of bounded linear operators on a
separable infinite-dimensional Hilbert space $\H$ and consider the ideal
 $$
{\cal L}^{1^+}:=
\left\{a\in B(\H):a{\rm\ is\ compact\ and\ }
\left({\s_n(a)\over\log n}\right)_{n\in\Na}
{\rm\ is\ a\ bounded\ sequence}\right\}
 $$
where $\s_n(a)=\sum_{k=0}^{n-1}\m_k(a)$ and $\{\m_k(a)\}_k$ is the sequence of
singular values of $a$ counted with multiplicities and arranged in a decreasing
order. In the following an operator $a$ is said $d^+$-summable (or
$a\in{\cal L}^{d^+}$) $d>0$, if $|a|^d\in{\cal L}^{1^+}$.

\begdefinition{1.1} A $K$-cycle is a triple $(\A, D,\H)$ where $\H$ a separable
Hilbert space, $\A$ is a $^*$-subalgebra of $B(\H)$, $D$ is an unbounded
selfadjoint operator on $\H$ such that $D$ has compact resolvent and $[D,a]$ is
bounded for every $a\in\A$.
 \np
A $K$-cycle is called $d^+$-summable if $|D|^{-1}$ is $d^+$-summable.
\enddefinition

\begremark{} In the previous definition $[D, a]$ ``bounded" means that one of
the following three equivalent conditions is satisfied (see e.g. [\rfr(BrRo1)],
proposition 3.2.55):
 \item{$(i)$} there exists a core $\D$ for $D$ (indeed the whole domain of $D$)
such that, for each $a\in\A$, the sesquilinear form
 $$
q(x,y):=(Dx,ay)-(a^*x,Dy)
 $$
is bounded on $\D\times\D$.
 \item{$(ii)$} $\A$ is contained in the domain of the derivation $i[D,\cdot]$,
which is, by definition, the infinitesimal generator of the one-parameter group
of automorphisms implemented by $e^{itD}$.
 \item{$(iii)$} For any $a\in\A$, $a\xi\in D(D)$ for any $\xi$ in
the domain $\D(D)$ of $D$ and $[D,a]$ is norm bounded on $\D(D)$.
\endremark
An important object in non commutative geometry is a particular
type of non normal traces on $B(\H)$ which where introduced by Dixmier in
[\rfr(Dixm1)] (see
also [\rfr(AGPS4),\rfr(GuIs1)] for further developments).
These traces are parametrized by two objects, one is the domain of the trace,
which is described in terms of the asymptotic behavior of the sequence
$\s_n(\cdot)$, the second is  a generalized limit procedure, more precisely a
dilation invariant state on bounded sequences on $\Na$. By {\it Dixmier trace}
is
generally meant a trace which ``sums'' logarithmic divergences, i.e. whose
domain is the ideal ${\cal L}^{1^+}$.

 We shall denote by $\t_\o$ the Dixmier trace corresponding to the state $\o$
(we shall not choose any particular state $\o$, but cf. also
[\rfr(Conn4)] about a possible canonical choice).

\bigskip

 \begtheorem{1.2} Let $(\A,D,\H)$ be a $d^+$-summable $K$-cycle. Then the
functional
 $$
\f(a):=\t_\o (|D|^{-d}a),\qquad \forall a\in B(\H)
 $$
is a hypertrace on $\A$, i.e.
 $$
\f(ab)=\f(ba),\qquad\forall a\in \A, b\in\B(\H)\ .
 $$
\endtheorem

Our proof relies on the following lemmas.

\beglemma{1.3} H\"older inequality holds for $\t_\w$ on $B(\H)$.
\endlemma
\begProof
Let $p,q>0$, $1/p+1/q=1$. It is sufficient to consider $a,b$ compact operators.
We have
 $$
\eqalign{\s_N(ab)
&=\sum_{k=0}^{N-1}\m_k(ab)\leq\sum_{k=0}^{N-1}\m_k(a)\m_k(b)\cr
&\leq\left(\sum_{k=0}^{N-1}\m_k(a)^p\right)^{1/p}
\left(\sum_{k=0}^{N-1}\m_k(b)^q\right)^{1/q}=
\left(\s_{N}(|a|^p)\right)^{1/p}\left(\s_{N}(|b|^q)\right)^{1/q}\cr}
 $$
by the Weyl inequality, H\"older inequality for $\Re^N$ and the spectral
theorem. Then, dividing by $\log N$ and applying the state $\w$ to the previous
inequality, we get, by definition of $\t_\w$,
 $$
\eqalign{\t_\w(|ab|)
&=\w\left({\s_N(ab)\over\log N}\right)\cr
&\leq\w\left(\left({\s_N(|a|^p)\over\log N}\right)^{1/p}
\left({\s_N(|b|^q)\over\log N}\right)^{1/q}\right)\cr
&\leq\w\left({\s_N(|a|^p)\over\log N}\right)^{1/p}
\w\left({\s_N(|b|^q)\over\log N}\right)^{1/q}
=\t_\w(|a|^p)^{1/p}\t_\w(|b|^q)^{1/q}\cr}
 $$
where we used H\"older inequality for states on abelian C$^*-$algebras.

The case $p=1$, $q=\infty$, $a$ compact and $b$ bounded can be proven by the
same methods.\endProof

The following lemma is well known, see e.g. [\rfr(CoMo1)]
[\rfr(Voic1)]. For the sake of completeness, we give a proof of the statement.

\beglemma{1.4} Let $a$ be a bounded operator and $D$ a selfadjoint operator
with
bounded inverse on $\H$ such that $[D,a]$ is bounded. Then, for any $0<r<1$,
the
operator $[|D|^r,a]$ is bounded and the following inequality holds:
 $$
\|[|D|^r,a]\|\leq C\|[D,a]\|.
 $$
where $C>0$ does not depend on $a$ and $\|\cdot\|$ denotes the usual operator
norm on $B(\H)$. \endlemma

\begProof
 First we show that the following property on the domain of a derivation holds
(cf. [\rfr(Brat1)]): let $D$ an unbounded selfadjoint
operator with bounded inverse, let $a$ be in the domain of the derivation
$[D,\cdot]$ and let $g$ be a $C^1$ function on $\Re$ such that $\widehat{g'}$,
the Fourier transform of $g'$, is in $L^1(\Re,dx)$. Then $a$ is in the domain
of
the derivation $[g(D),\cdot]$ and
 $$
\|[g(D),a\|\leq \|\widehat{g'}\|_1\|[D,a]\|.\eqno(1.1)
 $$
Let $\D$ be the space of vectors in $\H$ having compact support with respect
to the spectral measure of $D$. $\D$ is a common core for $D$ and $g(D)$ and,
for each $x,y\in\D$, the following formula holds (cf. again [\rfr(Brat1)]):
 $$
\eqalign{
&(g(D)x,ay)-(x,ag(D)y)=\cr
&={i\over\sqrt{2\pi}}\int_{-\infty}^{+\infty}dp\ p\hat g(p)\int_0^1 dt\
\left((De^{-itpD}x,ae^{i(1-t)pD}y)-(x,e^{itpD}aDe^{i(1-t)pD}y)\right)\cr}
 $$
 By a straightforward estimate on the previous equality and  the equivalence
stated in the remark following definition~1.1 we get that $[g(D),a]$ is bounded
and (1.1) holds.

 Finally let us choose $g\in C^\infty(\Re)$ such that $g(x)=|x|^r$ when
$|x|\geq dist(0,\s(D))$. For such a $g$ one has $g(D)=|D|^r$.

The lemma is proved if we show that the Fourier transform of $g'$ is in
$L^1(\Re,dx)$.

 Since $g'$ is the sum of an homogeneous function of degree $r-1$ plus a
function with compact support, its Fourier transform is the sum of an
homogeneous function of degree $-r$ plus an analytic function, therefore is
locally summable since $r<1$. But $g'$ is also $C^\infty$, therefore its
Fourier
transform is a function of rapid decay.\endProof

\begProofof{ Theorem 1.2} First let us notice that, making use of the trace
property of $\t_\o$, the statement of the theorem is equivalent to prove
 $$
 \t_\w\left(\left|[|D|^{-d},a]\right|\right)=0\qquad \forall a\in\A.\eqno(1.2)
 $$
 Then we recall that, if $H$ is a self-adjoint operator with bounded inverse
and $[H,a]$ is bounded, for any $k\in\Na$ the following identity holds:
 $$
[a,H^{-k}]=\sum_{j=1}^{k} H^{-j}[H,a]H^{-k-1+j}\eqno(1.3)
 $$
Now choose $r\in(0,1)$ such that $k:={d\over r}$ is a natural number. Then,
applying identity (1.3) with $H=|D|^r$, H\"older inequality and Lemma~1.4 we
get
 $$
\eqalign{
\t_\w\left(|[a,|D|^{-d}]|\right)&=\t_\w\left(|[a,(|D|^r)^{-k}]|\right)\cr
&\leq\sum_{i=1}^{k}
\t_\w\left(\left||D|^{-rj}[a,|D|^r]|D|^{r(-k-1+j)}\right|\right)\cr
&\leq\|[a,|D|^r]\|\sum_{j=1}^k\left(\t_\w\left(|D|^{-
rjp_j}\right)\right)^{1/p_j}
\left(\t_\w\left(|D|^{r(-k-1+j)q_j}\right)\right)^{1/q_j}\cr
&\leq C\|[a,D]\|\sum_{j=1}^k\left(\t_\w\left(|D|^{-rjp_j}\right)\right)^{1/p_j}
\left(\t_\w\left(|D|^{r(-k-1+j)q_j}\right)\right)^{1/q_j}\cr
}
 $$
 With a suitable choice of $p_j$ and $q_j$, e.g. $p_j={2d\over
r(2j-1)}$, $q_j={2d\over r(2k+1-2j)}$, both the exponents of $|D|$ in the last
term are strictly smaller than $-d$, therefore the Dixmier traces vanish, which
proves the theorem.\endProof

We refer to [\rfr(Conn4)] for conditions on the non triviality of $\f$ on $\A$.

As we pointed out in the introduction, one is interested in proving that the
state $\f$ constructed above gives a trace not only on the algebra
$\A$,
but also on the universal graded differential algebra $\Q(\A)$ on $\A$ via the
formula
 $$
\t(a_0da_1\dots da_n):=i^n\f(a_0[D,a_1]\dots[D,a_n])\eqno(1.4)
 $$
 This result does not hold in  full generality, as it is shown by the
counterexample described in the next section. Here we discuss regularity
conditions which ensure that the trace property holds.

As shown in [\rfr(Conn4)], the the condition $[D,\A]$ bounded
contained in the $K$-cycle assumption corresponds to the Lipshitz
regularity for functions.
According to [\rfr(Conn4), p.546], higher regularity conditions
may be given in terms of the derivation $\d$ corresponding to the
commutator with $|D|$, namely the generator of the following
one-parameter group of automorphisms of $B(\H)$:
 $$
 \a_t(a):=e^{it|D|}ae^{-it|D|}\ ,\qquad a\in B(\H).\eqno(1.5)
 $$
 More precisely we may introduce the subalgebras $\A_n$ generated
by the elements $a\in\A$ such that $a$ and $[D,a]$ are in the
domain of $\d^{n-1}$. The intersection of such algebras corresponds
to $\C^\infty$ functions.

 First we notice that while Lipshitz regularity is expressed by the
boundedness of the commutator with $D$, higher commutators with
$D$ do not describe higher regularity, because in the commutative
setting $[D,[D,f]]$ is a sum of a Clifford multiplication operator
and a differential operator, and therefore cannot be bounded.
 On the other hand, $\d^n$ is bounded on $\C^\infty$ functions (or
forms) in the commutative case.

This fact may be explained since the flow described in (1.5) is
deeply related to the geodesic flow (cf. [\rfr(DuGu1)],
[\rfr(Bess1)]).

Now we show that the mentioned regularity conditions imply the trace property
on
$\Q(\A)$. First we remark that the functional $\t$ on $\Q(\A)$ in formula (1.4)
is a trace if and only if $\f$ is a trace on the $^*$-algebra
$\langle\A,[D,\A]\rangle$ generated by $\A$ and the commutators of
$D$ with the elements of $\A$. As a consequence the following
property holds.

\begproposition{1.5} The functional $\t$ is a trace on $\Q(\A_2)$.
\endproposition

\begProof In fact $(\langle\A_2,[D,\A_2]\rangle,|D|,\H)$ is a
$d^+$-summable $K$-cycle, hence the assert follows from theorem~1.2
and the previous observation.\endProof

We may restate the proposition saying that $\t$ is a trace on
$\Q(\A)$ if, extending the usual notation for Sobolev spaces, the
$K$-cycle is $W^{2,\infty}$ regular, namely $\A=\A_2$.

However, the result given by Proposition~1.5 is
not completely analogous to its commutative counterpart. In fact,
the proof of the trace property in the  commutative case (cf. e.g.
[\rfr(GBVa1)]) is based on the fact that, since the principal
symbol of the product of two pseudo-differential operators is the
product of the principal symbols of the factors, the order of the
commutator is the sum of the orders of the factors minus 1.
Therefore, on an $n$-dimensional manifold, the Dixmier trace
vanishes on the commutator $[|D|^{-n},\s]$, when $\s$ is a section
of the Clifford bundle, because it is a pseudo-differential
operator of order $n-1$, and the regularity (differentiability)
assumptions on $\s$ plays no crucial role.

A natural assumption, that hold in the commutative case and allows
to extend the trace property to $\Q(\A)$, is to ask that $\A_2$ is
a large enough subalgebra of $\A$.

 \begtheorem{1.6} Let $(\A,D,\H)$ be a $d^+$-summable $K$-cycle.
If $\A_2$ is a core for the derivation $[D,\cdot]$ on $\A$, then $\t$ given by
(1.4) is a trace on the algebra $\Q(\A)$.
\endtheorem

	\begProof As already remarked it is sufficient to show the
trace property of $\f$ on $\langle\A,[D,\A]\rangle$. By
Proposition~1.5 $\f$ is a trace on $\langle\A_2,[D,\A_2]\rangle$.
Since $\A_2$ is dense in $\A$ in the graph norm of $[D,\cdot]$
(hence in the usual norm, too), $[D,\A_2]$ is norm dense in
$[D,\A]$ and this implies the trace property on (the closure of)
$\langle\A,[D,\A]\rangle$.\endProof

\titlea{2. An example of a K-cycle}

Let $\K$ be an infinite dimensional Hilbert space and $b$ an unbounded
selfadjoint positive operator with compact inverse, such that
$0<\tau_{\omega}(b^{-d})<\infty$. We consider on $\H=\Co^2\otimes \K$
the $^*$-algebra $\A$ generated by  the element $a$ in $B(\H)$, where
 $$
a=m_{12}\otimes b^{-1},
 $$
and $m_{ij}$, $i,j=1,2$, are the matrix units.

We consider also the unbounded selfadjoint operator $D\equiv \alpha\otimes b$
on $\H$, where $\alpha=\left(\matrix{\lambda &0\cr0&\mu}\right)$, $\lambda$,
$\mu\in\Re$.

\begproposition{2.1} The triple $(\A,D,\H)$ is a $d^+$-summable $K$-cycle.
\endproposition

\begProof
Since $a^2=(a^*)^2=0$, $\A$ is linearly generated by the following elements:
 $$
x_{11}^k\equiv(a^*a)^{k+1},\quad x_{22}^k\equiv (aa^*)^{k+1},\quad
x_{12}^k\equiv a(a^*a)^k,\quad x_{21}^k\equiv a^*(aa^*)^k,\qquad k\ge 0.
 $$
A computation shows:
 $$
\eqalign{
x_{ii}^k&=m_{ii}\otimes b^{-(2k+2)}\qquad i=1,2 \cr
x_{ij}^k&=m_{ij}\otimes b^{-(2k+1)}\qquad i\ne j.\cr}
 $$
 Denoting by $Tr_{\o}$ the Dixmier trace on $\H$, we have
$Tr_\o(|D|^{-d})=tr(|\a|^{-d})\cdot\t_{\o}(b^{-d})<\infty$
and moreover
 $$
[D,x_{ij}^k]=(i-j)(\l-\m) m_{ij}\otimes
b^{-2k}\qquad i,j=1,2
 $$
hence $[D,\A]\subseteq B(\H)$.\endProof

Theorem~1.2 implies that $\f(\cdot)\equiv Tr_{\o}(\cdot |D|^{-d})$ is a trace
on $\A$. The following proposition completely characterizes the trace
property on $\Q(\A)$ in terms of $\l$ and $\m$.

 \begproposition{2.2} The functional $\t$ on $\Q(\A)$ is a trace on
if and only if $|\l|=|\m|$.
 \endproposition

\begProof As explained in the preceding section, this is equivalent to show
that $\f$ is a trace on $\langle\A,[D,\A]\rangle$.

 If $|\l|=|\m|$, the state
$\f$ writes as $\f(\cdot)=|\l|^{-d}tr(\cdot)\otimes\t_\w(\cdot b^{-d})$. Since
$\langle\A,[D,\A]\rangle$ is contained in $M_2\otimes\C(\s(b^{-1}))$ the trace
property of $\f$ follows.

Conversely, let us assume that $\f$ is a trace on
$\langle\A,[D,\A]\rangle$.
Since we have
 $$
[[D,a],[D,a^{*}]]=-(\l-\m)^2(m_{11}-m_{22})\otimes1\ ,
 $$
it follows
 $$
0=\f([[D,a],[D,a^{*}]])=
-(\l-\m)^2 (|\l|^{-d}-|\m|^{-d})\t_\o(b^{-d})
 $$
which proves the assert.\endProof

It is instructive to verify by direct arguments that, when
$|\l|\ne|\m|$, $\A_2$ is not dense in $\A$ in the graph norm of the
derivation $[D,\cdot]$. More precisely, motivated by theorem~1.6,
we compute the norm closures of $\A$ and $\A_2$ and their closures
in the graph norm of the derivation $[D,\cdot]$, showing that while
the first two closures coincide, the latter do not.

First we notice that the maps
 $$
p\to m_{ii}\otimes p(b^{-1}),\qquad i=1,2
 $$
are isomorphisms from the algebra of even polynomials of degree $\geq2$ on
$\s(b^{-1})$ to the diagonal blocks of $\A$ (which coincide with the diagonal
blocks of $\A_2$). Therefore the diagonal blocks of the norm closure of
$\A$ coincide with the corresponding diagonal blocks of the norm closure of
$\A_2$ and are isomorphic to the norm closure of the algebra of even
polynomials
of degree $\geq2$ on $\s(b^{-1})$ which is $\C_0(\s(b^{-1}))$, the algebra of
continuous functions on the spectrum of  $b^{-1}$ vanishing at 0.

In a similar way the maps
 $$
p\to m_{ij}\otimes p(b^{-1}),\qquad i\ne j,\  i,j=1,2
 $$
give isomorphisms from the vector space of odd polynomials of degree $\geq1$ on
$\s(b^{-1})$ to the anti-diagonal blocks of $\A$ and isomorphisms from the
vector space of odd polynomials of degree $\geq3$ on $\s(b^{-1})$ to the
anti-diagonal blocks of $\A_2$. Since the norm closure of the vector space of
odd polynomials of degree $\geq p>0$ on $\s(b^{-1})$ is $\C_0(\s(b^{-1}))$, the
anti-diagonal blocks of the norm closure of $\A$ coincide with the
corresponding
anti-diagonal blocks of the norm closure of $\A_2$ and are isomorphic to
$\C_0(\s(b^{-1}))$. Finally we get
$\ov\A=\ov\A_2\simeq M_2\otimes\C_0(\s(b^{-1}))$.

Now we consider the closures of $\A$ and $\A_2$ in the graph norm
of $[D,\cdot]$. Since the commutator with $D$ is zero on diagonal elements and
coincides up to a constant with the multiplication by $1\otimes b$ on the
anti-diagonal elements, the diagonal part of the closures of $\A$ and of $\A_2$
in the graph norm are equal to the corresponding norm closures and
therefore coincide and, since $b^{-1}$ is bounded, the graph norm on the
antidiagonal elements $m_{ij}\otimes p(b^{-1})$ is equivalent to the norm
 $$
\norm{m_{ij}\otimes p(b^{-1})}:=\|bp(b^{-1})\|.
 $$
 Then the maps
 $$
p\to m_{ij}\otimes b^{-1}p(b^{-1}),\qquad i\ne j,\  ij=1,2
 $$
gives isometric isomorphisms between the vector space of even polynomials on
the
spectrum of $b^{-1}$ and the anti-diagonal blocks of $(\A,\norm\cdot)$ and
isometric isomorphisms between the vector space of even polynomials of degree
$\geq2$ on the spectrum of $b^{-1}$ and the anti-diagonal blocks of
$(\A_2,\norm\cdot)$.

As a consequence, an anti-diagonal block of the graph norm closure
of $\A$ is isomorphic to the space of continuous functions on
$\s(b^{-1})$ vanishing at $0$ and with finite derivative in $0$,
while an anti-diagonal block of the graph norm closure of $\A_2$ is
isomorphic to the space of continuous functions on $\s(b^{-1})$
vanishing at $0$ and with null derivative in $0$.

\titlea{Acknowledgments} We thank Prof.A.Connes for very useful
remarks and suggestions on the first version of this work.

%%% end of text

\references
\end